\documentclass[a4paper,11pt]{article}
\usepackage{subcaption}
\usepackage{float}
\usepackage{epsfig}
\usepackage{amsmath}
\usepackage{amsfonts}
\usepackage{amsthm}
\usepackage{tikz}
\usetikzlibrary{arrows}

\newtheorem{theorem}{{\sc Theorem}}[section]
\newtheorem{lemma}{{\sc Lemma}}[section]

\title{Discontinuous Galerkin Finite Element Methods for 1D Rosenau Equation}

\author{
P. Danumjaya  \footnote{e-mail address:
  danu@goa.bits-pilani.ac.in} \hspace{0.5in} K. Balaje \\
  Department of Mathematics, \\
  BITS-Pilani K K Birla Goa Campus, \\
  Goa-403726, India.
}
\date{}

\begin{document}
\maketitle

\begin{abstract}
  \noindent
  In this paper, discontinuous Galerkin
  finite element methods are applied to one dimensional Rosenau
  equation. Theoretical results including 
  consistency, \textit{a priori}
  bounds and optimal error estimates are established for both
  semidiscrete and fully discrete schemes. Numerical
  experiments are performed to validate the theoretical   
  results. The decay estimates are verified numerically for the
  Rosenau equation. 
\end{abstract}

\section{Introduction}

Consider the one dimensional Rosenau equation:
\begin{equation}
  u_t + u_{xxxxt} = f(u)_x, \quad (x,t) \in (a,b) \times (0,T] \label{Eq1.1}
\end{equation}
with initial condition
\begin{equation}
  u(x,0) = u_0(x), \quad x \in (a, b), \label{Eq1.2}\\
\end{equation}
and the boundary conditions
\begin{eqnarray}
  u(a, t) &=& u(b, t) = 0,\nonumber \\
  u_x(a, t) &=& u_x(b, t) = 0, \label{Eq1.3}
\end{eqnarray}
where $f(u)$ is a nonlinear term in $u$ of the type $f(u) = \displaystyle
\sum_{i=1}^{n} \frac{c_i u^{p_i+1}}{p_i+1}$, here $c_i$ is a real
constant and $p_i$ is a positive integer. \\\\
The Rosenau equation is an example of a nonlinear partial differential
  equation, which governs
  the dynamics of dense discrete systems and models wave propagation in
  nonlinear dispersive media. \\\\
Recently, several numerical techniques
  like conforming finite element methods, mixed finite element methods,
  orthogonal cubic spline collocation methods, etc., were proposed to
  find the approximate solution of Rosenau equation.  
The different conforming finite element techniques are used to
  approximate the solution of Rosenau equation needs $C^1$-interelement
  continuity condition and mixed finite element formulations
  are required $C^0$-continuity condition. In this article discontinuous
Galerkin finite element methods are used to approximate the solution. \\\\
The well-posedness results of (\ref{Eq1.1})-(\ref{Eq1.3}) was
proved by Park \cite{park} and Atouani {\it et al.} in \cite{atou}. 
Earlier, some numerical methods were proposed to solve the Rosenau
equation (\ref{Eq1.1})-(\ref{Eq1.3}) using finite difference methods
by  Chung \cite{chung3}, conservative difference schemes by Hu and Zheng
\cite{Hu2010} and Atouni and Omrani \cite{omrani}.
Finite element Galerkin method was used by
\cite{atou,ha}, a second order splitting combined with orthogonal
cubic spline collocation method was used by Manickam {\it et al.} \cite{manickam} and Chung
and Pani in \cite{chung} constructed
a $C^1$-conforming finite element method for the Rosenau equation
(\ref{Eq1.1})-(\ref{Eq1.3}) in
two-space dimensions. \\\\
In recent years, there has been a growing interest in discontinuous
Galerkin finite element methods because of their flexibility in
approximating globally rough solutions and their potential for error
control and mesh adaptation.\\\\
Recently, a cGdG method was proposed by Choo. {\it et. al}
in \cite{dgros}. A subdomain finite element method using sextic
b-spline was proposed by Battal and Turgut in \cite{subd}. But constructing
$C^1$ finite elements for fourth order problems becomes expensive
and hence discontinuous Galerkin finite element methods can be used to
solve fourth order problems \cite{gudi}.\\\\
In this paper, we introduce discontinuous Galerkin finite element
methods (DGFEM) in space to solve the one dimensional Rosenau
equation (\ref{Eq1.1})-(\ref{Eq1.3}). Comparitive to existing methods
our proposed method require less regularity. \\\\
The outline of the paper is as follows. In Section 2, we derive the
discontinuous weak formulation of the Rosenau equation. In Section 3,
we discuss the \textit{a priori} bounds and optimal error estimates
for the semidiscrete problem. In Section 4, we discretize the
semidiscrete problem in the temporal direction using a backward Euler
method and discuss the \textit{a priori} bounds and optimal error
estimates. Finally, in Section 5, we present some numerical results to
validate the theoretical results. \\\\
Throughout this paper, $C$ denotes a generic positive constant which
is independent of the discretization parameter $h$ which may have
different values at different places.  
\section{Weak Formulation}

In this section, we derive the weak formulation for the problem
(\ref{Eq1.1})-(\ref{Eq1.3}). \\\\
We discretize the domain $(a, b)$ into $N$ subintervals as
$$
a = x_0 < x_1 < x_2 < \dots <x_N = b,
$$
and $I_n = \left(x_n, x_{n + 1} \right)$ for $n = 0, 1, 2, \ldots,
N-1$. 
We denote this partition by $\mathcal{E}_h$ consisting of sub-intervals
$I_n, \, n = 0, 1, 2, \ldots N-1$.  Below, we define
the broken Sobolev space and corresponding norm 
\begin{equation*}
  H^s(\mathcal{E}_h) = \left\{ v \in L^2(a, b) \;: \; v |_{I_n}
    \in H^s(I_n), \; \forall I_n\in\mathcal{E}_h \right\}
\end{equation*}
and
\begin{equation*}
  |||v|||_{H^s(\mathcal{E}_h)} = \left(\sum_{K\in\mathcal{E}_h} \lVert
    v \rVert_{H^s(K)}^2\right)^\frac{1}{2}.
\end{equation*}
\noindent
Now define the jump and
average of $v$ across the nodes $\{x_n\}_{n=1}^{N-1}$ as follows.
The jump of a function value $v(x_n)$ across
the inter-element node $x_n$, shared by $I_{n-1}$ and $I_n$ denoted
by $[v(x_n)]$ and defined by
\begin{equation*}
  [v(x_n)] = v(x_n^-) - v(x_n^+).
\end{equation*}
At the boundary $x_0$ and $x_N$, we set
\begin{equation*}
  [v(x_0)] = -v(x_0), \;\;\;\; \mbox{and} \;\;\;\;
  [v(x_N)] = v(x_N).
\end{equation*}
\noindent
The average of a function value $v(x_n)$ across
the inter-element node $x_n$, shared by $I_{n-1}$ and $I_n$ denoted
by $\left\{ v(x_n) \right\}$ and defined by
\begin{equation*}
  \{v(x_n)\} = \frac{1}{2} \left(v(x_n^-) + v(x_n^+)\right).
\end{equation*}
At the boundary $x_0$ and $x_N$, we set
$$
\{v(x_0)\} = v(x_0),  \;\;\;\; \mbox{and} \;\;\;\;
\{v(x_N)\} = v(x_N).
$$
\noindent
We multiply \eqref{Eq1.1} with $v\in H^s(\mathcal{E}_h)$ and
integrate over $I_n = (x_n,x_{n+1})$ to obtain
\begin{eqnarray}
  \int_{x_n}^{x_{n+1}} (u_t + u_{xxxxt}) v \; dx =
  \int_{x_n}^{x_{n+1}} f(u)_x v\; dx. \label{Eq2.1}
\end{eqnarray}
Now, using integration by parts twice in (\ref{Eq2.1}), we arrive at
\begin{eqnarray*}
  \int_{x_n}^{x_{n+1}} u_t v \; dx  + \int_{x_n}^{x_{n+1}} u_{xxt} v_{xx} \; dx
  + u_{xxxt}(x_{n+1})v(x_{n+1}^-) - u_{xxxt}(x_{n})v(x_{n}^+) \\ - u_{xxt}(x_{n+1})v_x(x_{n+1}^-)
  + u_{xxt}(x_{n})v_x(x_{n}^+) = \int_{x_n}^{x_{n+1} } f(u)_x v \; dx.
\end{eqnarray*}
Summing over all the elements $n = 0,1,\dots,N-1$ and using
$$
ps-qr = \frac{1}{2}(p+q)(s-r) + \frac{1}{2}(r+s)(p-q), \;\;\; p, q,
r \; \mbox{and} \; s \in \mathbb{R},
$$
we obtain
\begin{eqnarray}
  \sum_{n=0}^{N-1} \int_{x_n}^{x_{n+1}} u_t v \; dx
  + \sum_{n=0}^{N-1} \int_{x_n}^{x_{n+1}} u_{xxt} v_{xx} \; dx
  &+& \sum_{n=0}^{N} \big\{u_{xxxt}(x_n)\big\}\big[v(x_n)\big] -
      \sum_{n=0}^{N} \big\{u_{xxt}(x_n)\big\}\big[v_x(x_n)\big]
      \nonumber \\
  &=& \sum_{n=0}^{N-1} \int_{x_n}^{x_{n+1}} f(u)_x v \; dx. \label{Eq2.2}
\end{eqnarray}
Since $u(x,t)$ is assumed to be sufficiently smooth, we have
$\big[u_t(x_n)\big] = \big[u_{xt}(x_n)\big] = 0$. Using this, we write as
\begin{multline}
  \sum_{n=0}^{N}\big\{v_{xxx}(x_n) \big\}\big[u_t(x_n)\big] -
  \sum_{n=0}^N\big\{v_{xx}(x_n)\big\}\big[u_{xt}(x_n)\big]  +
  \sum_{n=0}^N\frac{\sigma_0}{h^\beta}\big[v(x_n)\big]\big[u_t(x_n)\big] +
  \sum_{n=0}^N\frac{\sigma_1}{h}\big[v_x(x_n)\big]\big[u_{xt}(x_n)\big] \\
  = v_{xxx}(x_0)u_t(x_0) - v_{xxx}(x_N)u_t(x_N) -
  v_{xx}(x_0)u_{xt}(x_0) + v_{xx}(x_N)u_{xt}(x_N) -\\
  \frac{\sigma_0}{h^\beta}v(x_0)u_t(x_0) -
  \frac{\sigma_0}{h^\beta}v(x_N)u_{t}(x_N) -
  \frac{\sigma_1}{h}v_{x}(x_0)u_{xt}(x_0) -
  \frac{\sigma_1}{h}v_{x}(x_N)u_{xt}(x_N) = 0. \label{Eq2.3}
\end{multline}
The right hand side of \eqref{Eq2.3} was found out using the
boundary conditions \ref{Eq1.3}. Adding \eqref{Eq2.3} to \eqref{Eq2.2} we obtain
\begin{eqnarray}
  \sum_{n=0}^{N-1} \int_{x_n}^{x_{n+1}} u_t v \; dx
  &+& \sum_{n=0}^{N-1} \int_{x_n}^{x_{n+1}} u_{xxt} v_{xx} \; dx
      + \sum_{n=0}^{N} \big\{u_{xxxt}(x_n)\big\}\big[v(x_n)\big]  +
      \sum_{n=0}^{N}
      \big\{v_{xxx}(x_n)\big\}\big[u_t(x_n)\big]\nonumber \\
  &-& \sum_{n=0}^{N} \big\{u_{xxt}(x_n)\big\}\big[v_x(x_n)\big]
      - \sum_{n=0}^{N}
      \big\{v_{xx}(x_n)\big\}\big[u_{xt}(x_n)\big]
      + \frac{\sigma_0}{h^\beta} \sum_{n=0}^{N}
      \big[v(x_n)\big]\big[u_t(x_n)\big] \nonumber \\
  &+& \frac{\sigma_1}{h} \sum_{n=0}^{N} \big[v_x(x_n)\big]\big[u_{xt}(x_n)\big]
      = \sum_{n=0}^{N} \int_{x_n}^{x_{n+1}} f(u)_x v \; dx \label{Eq2.4}.
\end{eqnarray}
We define the bilinear form as
\begin{eqnarray}
  B(u,v) &=& A(u,v) + J^{\sigma_0}(u,v) + J^{\sigma_1}(u,v),\label{Eq:n2.4}
\end{eqnarray}
where
\begin{eqnarray*}
  A(u,v) &=& \sum_{n=0}^{N-1}\int_{x_n}^{x_{n+1}} u_{xx}v_{xx} \,dx +
             \sum_{n=0}^N\left(\big\{ u_{xxx}(x_n)\big\}\big[v(x_n)\big] + \big\{
             v_{xxx}(x_n)\big\}\big[u(x_n)\big]\right)\\
         &-&\sum_{n=0}^N\left(\big\{ u_{xx}(x_n)\big\}\big[v_x(x_n)\big] + \big\{
             v_{xx}(x_n)\big\}\big[u_x(x_n)\big]\right),
\end{eqnarray*}
and
\begin{eqnarray*}
  J^{\sigma_0}(u,v) =
  \sum_{n=0}^N\frac{\sigma_0}{h^\beta}\big[u(x_n)\big]\big[v(x_n)\big],
  \qquad \qquad J^{\sigma_1}(u,v) =
  \sum_{n=0}^N\frac{\sigma_1}{h}\big[u_x(x_n)\big]\big[v_x(x_n)\big].
\end{eqnarray*}
In \eqref{Eq:n2.4}, $J^{\sigma_0}$ and $J^{\sigma_1}$ are the penalty terms
and $\sigma_0, \sigma_1 > 0$. The value of $\beta$ will be defined
later. \\\\
The weak formulation of
\eqref{Eq1.1}-\eqref{Eq1.3} as follows:
Find $u(t) \in H^s(\mathcal{E}_h), \; s > 7/2$, such that
\begin{eqnarray}
  (u_t,v) + B(u_t,v) &=& \left( f(u)_x,v\right), \; \forall v \in
                         H^s(\mathcal{E}_h), \;\;t>0 \label{Eq2.7} \\
  u(x,0) &=& u_0(x). \label{Eq2.8}
\end{eqnarray}
\noindent
Below, we state and prove the consistency result of the weak formulation
(\ref{Eq2.7})-(\ref{Eq2.8}).

\begin{theorem}
  Let $u(x,t)\in C^4(a,b)$ be a solution of the continuous
  problem \eqref{Eq1.1}-\eqref{Eq1.3}. Then $u(x,t)$ satisfies the
  weak formulation \eqref{Eq2.7}-\eqref{Eq2.8}. Conversely, if
  $u(x,t)\in H^2(a, b) \cap H^s(\mathcal{E}_h)$ for $s>7/2$ is a
  solution of \eqref{Eq2.7}-\eqref{Eq2.8}, then $u(x, t)$ satisfies
  \eqref{Eq1.1}-\eqref{Eq1.3}.
  \begin{proof}
    Let $u(x,t)\in C^4(a,b)$ and $v\in
    H^s(\mathcal{E}_h)$. Multiply \eqref{Eq1.1} by $v$ and integrate
    from $x_n$ to $x_{n+1}$. Sum over all
    $n=0,1,\dots,N-1$ and using \eqref{Eq2.2}, \eqref{Eq2.3} and
    \eqref{Eq2.4}, we obtain the weak formulation \eqref{Eq2.7}.\\\\
    Conversely, let $u\in H^2(a,b) \cap H^s(\mathcal{E}_h), \; s>7/2$ and $v
    \in \mathcal{D}(I_n)$, the space of infinitely differentiable
    functions with compact support in $I_n$. Then, \eqref{Eq2.7}
    becomes
    \begin{equation}
      \int_{x_n}^{x_{n+1}} u_t v \; dx + \int_{x_n}^{x_{n+1}} u_{xxt}
      v_{xx} \; dx = \int_{x_n}^{x_{n+1}} f(u)_x v\; dx. \label{eq:Neq21}
    \end{equation}
    Applying integration by parts twice on the second term on the left hand
    side of \eqref{eq:Neq21} to obtain,
    \begin{equation*}
      \int_{x_n}^{x_{n+1}} u_{xxt} v_{xx} \; dx = \int_{x_n}^{x_{n+1}}
      u_{xxxxt} v \; dx,
    \end{equation*}
    as $v$ is compactly supported on $I_n$. This immediately yields
    \begin{equation}
      u_t + u_{xxxxt} = f(u)_x, \quad \text{a.e in}  \;\;I_n.\label{eq:3-17}
    \end{equation}
    Consider the node $x_k$ shared between $I_{k-1}$ and $I_k$. Choose
    $v \in H^2_0(I_{k-1}\cup I_{k})$, multiply \eqref{eq:3-17} by $v$
    and integrate over $(a,b)$ to obtain
    \begin{equation}
      \int_{I_{k-1} \cup I_k} u_t v \; dx + \int_{I_{k-1} \cup I_k}
      u_{xxxxt} v \; dx = \int_{I_{k-1} \cup I_k} f(u)_x v\;
      dx.\label{eq:Eq22}
    \end{equation}
Applying integration by parts twice on the second term of
    \eqref{eq:Eq22} and  using $v\in H^2_0(I_{k-1}\cup I_{k})$, we obtain
    \begin{eqnarray}
      \int_{I_{k-1} \cup I_k} u_t v \; dx
      &+& \int_{I_{k-1} \cup I_k}
          u_{xxt} v_{xx} \; dx + \left[u_{xxxt}(x_k)\right]v(x_k)
          -
          \left[u_{xxt}(x_k)\right]v_x(x_k)
          \nonumber \\
      &=&  \int_{I_{k-1} \cup I_k} f(u)_x v \, dx. \label{eq:3-18}
    \end{eqnarray}
    On the other hand, we have from \eqref{Eq2.7} for the choice of $u$ and $v$,
    \begin{eqnarray}
      \int_{I_{k-1} \cup I_k} u_t v \; dx + \int_{I_{k-1} \cup I_k}
      u_{xxt} v_{xx} \; dx = \int_{I_{k-1} \cup I_k} f(u)_x v\, dx. \label{eq:3-19}
    \end{eqnarray}
    Comparing \eqref{eq:3-18} and \eqref{eq:3-19} and using the fact
    that $v$ is arbitrary, we obtain
    \[\left[u_{xxxt}(x_k)\right] = 0.\]
    Thus $u_{xxxxt} \in L^2(\Omega)$ and hence, from \eqref{eq:Eq22},
    we obtain
    \begin{equation}
      u_t + u_{xxxxt} = f(u)_x \qquad \text{a.e in  } (a,b).
    \end{equation}
    This completes the proof.
  \end{proof}

\end{theorem}

\section{Semidiscrete DGFEM}
\setcounter{equation}{0}
In this section, we discuss the {\it a priori} bounds and optimal
error estimates for the semidiscrete Galerkin method. \\\\
We define a finite dimensional subspace $D_k(\mathcal{E}_h)$ of
$H^s(\mathcal{E}_h), \; s > 7/2$ as
\begin{equation*}
  D_k(\mathcal{E}_h) = \left\{ v\in L^2(a, b): v |_{I_n}
    \in \mathbb{P}_k(I_n), \;\;I_n \in \mathcal{E}_h\right\}.
\end{equation*}
The weak formulation for the semidiscrete Galerkin method is to find
$u^h(t) \in D_k(\mathcal{E}_h)$ such that
\begin{eqnarray}
  (u^h_t,\chi) + B(u^h_t,\chi) &=& \left( f(u^h)_x,\chi\right), \; \; \mbox{for
                                   all} \;\;  \chi \in D_k(\mathcal{E}_h), \label{Neq2.1}\\
  u^h(0) &=& u^h_0, \label{Neq2.2}
\end{eqnarray}
where $u_0^h$ is an appropriate approximation of $u_0$ which will be
defined later.

\subsection{\textit{A priori} Bounds}
In this sub-section, we derive the \textit{a priori} bounds. \\

\noindent
Define the energy norm
$$
||u||_{\mathcal{E}}^2 = \sum_{n=0}^N\int_{x_n}^{x_{n+1}} u_{xx}^2 + \sum_{n=0}^N
\frac{\sigma_0}{h^\beta} |[u(x_n)]|^2 + \sum_{n=0}^N \frac{\sigma_1}{h}
|[u_x(x_n)]|^2.
$$
We note from \cite{gudi} that $B(.,.)$ is coercive with respect to the energy
norm, i.e.,
$$
B(v,v) \ge C||v||_{\mathcal{E}}^2, \quad v\in D^k(\mathcal{E}_h),
$$
for sufficiently large values of $\sigma_0$ and $\sigma_1$.\\\\
Observe that \eqref{Neq2.1} yields a system of non-linear ordinary
differential equations and the existence and uniqueness of the solution can be
guaranteed locally using the Picard's theorem. To obtain existence and
uniqueness globally, we use continuation arguments and hence we need the
following \textit{a priori} bounds.

\begin{theorem}
  Let $u^h(t)$ be a solution to \eqref{Neq2.1} and assume that $f'$ is
  bounded. Then there exists a positive constant $C$ such that
  \begin{equation}
    \|u^h(t)\| + \|u^h(t)\|_{\mathcal{E}} \le C(\|u_0^h\|_{\mathcal{E}}).
  \end{equation}

  \begin{proof}
    On setting $\chi = u^h$ in \eqref{Neq2.1}, we obtain
    \begin{equation}
      (u^h_t,u^h) + B(u^h_t,u^h) = (f(u^h)_x,u^h).\label{Neq2.9}
    \end{equation}
    We rewrite the equation \eqref{Neq2.9} as
    $$
    \frac{1}{2}\frac{d}{dt}\|u^h(t)\|^2 + \frac{1}{2}\frac{d}{dt}
    B(u^h(t),u^h(t)) = (f'(u^h)u^h_x,u^h).
    $$
    Integrating from $0$ to $t$, we obtain
    \begin{equation}
      \|u^h(t)\|^2 + B(u^h(t),u^h(t)) = \|u^h(0)\|^2 + B(u^h(0),u^h(0)) +
      \int_{0}^{t}(f'(u^h)u^h_x,u^h) \,ds \label{Eq:2}.
    \end{equation}
    On using the coercivity of $B(u^h(t),u^h(t))$ and the
    boundedness of $f'$, we arrive at
    \begin{equation}
      \|u^h(t)\|^2 + C\|u^h(t)\|_{\mathcal{E}}^2 \le \|u^h(0)\|^2 + B(u^h(0),u^h(0)) +
      C \int_{0}^{t}(u^h_x,u^h) \,ds. \label{Eq:4}
    \end{equation}
    Using the Cauchy Schwarz inequality and the
    Poincar\'e inequality on the right hand side of \eqref{Eq:4}, we
    obtain
    \begin{equation}
      \|u^h(t)\|^2 + C\|u^h(t)\|_{\mathcal{E}}^2 \le C(\|u^h_0\|_{\mathcal{E}}) +
      \int_{0}^{t}\|u^h(s)\|_{\mathcal{E}}^2 \,ds \label{Eq:5}.
    \end{equation}
    An application of Gronwall's inequality yields the desired \textit{a
      priori} bound for $u^h(t)$.
  \end{proof}
\end{theorem}

\subsection{Error Estimates in the energy and $L^2$-norm}
In this subsection, we derive the optimal error estimates in energy
and $L^2$-norm.\\\\
Often a direct comparison between $u$ and $u^h$ does
not yield optimal rate of convergence. Therefore, there is a need to introduce an
appropriate auxiliary or intermediate function $\tilde{u}$ so that the
optimal estimate of $u-\tilde{u}$ is easy to obtain and the
comparision between $u^h$ and $\tilde{u}$ yields a sharper estimate
which leads to optimal rate of convergence for $u-u^h$. In literature, Wheeler \cite{wheeler}
for the first time introduced this technique in the context of
parabolic problem. Following Wheeler \cite{wheeler}, we introduce $\tilde{u}$ be an
auxiliary projection of $u$ defined by
\begin{equation}
  B(u-\tilde{u},\chi) = 0, \;\; \mbox{for all} \; \chi \in
  D^k(\mathcal{E}_h).\label{Eq:10}
\end{equation}
Now set the error $e = u-u^h$ and split as follows: 
$e = u - \tilde{u} - \left(u^h -\tilde{u} \right) = \eta - \theta$,
where $\eta = u - \tilde{u}$ and $\theta = u^h-\tilde{u}$.
Below, we state some error estimates for $\eta = u-\tilde{u}$ and its temporal
derivative.
\begin{lemma}\label{Lemma3}
  For $t\in(0,T]$ and $s > 7/2$ then there exists a positive constant $C$
    independent of $h$ such that the following error estimates for $\eta$ hold:
  \begin{eqnarray*}
    \left\|\frac{\partial^l \eta}{\partial t^l}\right\|_{\mathcal{E}}
    &\le&
          Ch^{\min(k+1,s)-2}\left(\sum_{m=0}^{l} |||\frac{\partial ^m
          u}{\partial t^m}|||_{H^s(\mathcal{E}_h)}\right),\\
    \left\|\frac{\partial^l \eta}{\partial t^l}\right\|
    &\le&
          Ch^{\min(k+1,s)}\left(\sum_{m=0}^{l} |||\frac{\partial ^m
          u}{\partial t^m}||| \right), \quad l = 0,1.
  \end{eqnarray*}
  \begin{proof}
    We split $\eta$ as follows:
    $$
    \eta = u - \tilde{u} = \left(u - \bar{u} \right) - \left(
    \tilde{u} - \bar{u} \right) = \rho - \xi,
    $$
    where $\rho = u - \bar{u}$, $\xi = \tilde{u} - \bar{u}$ and
    $\bar{u}$ is an interpolant of $u$ satisfying good approximation
    properties. Now from \eqref{Eq:10}, we have
    \begin{equation}
      B(\xi,\chi) = B(\rho,\chi)\label{eq:10}.
    \end{equation}
    We note that $\rho$ satisfies the following approximation
    property \cite{rvg}:
    $$
    \|\rho\|_{H^q(I_n)} \le Ch^{\min(k+1,s)-q} \|u\|_{H^s(I_n)}.
    $$
    Set $\chi = \xi$ in
    \eqref{eq:10} to obtain
    $$
    B(\xi,\xi) = B(\rho,\xi).
    $$
    A use of coercivity of $B(.,.)$ and the assumption that $\bar{u}$ is a
    sufficiently smooth interpolant of $u$, we obtain
    \begin{equation}
      C \|\xi\|^2_{\mathcal{E}} \le \sum_{n=0}^{N-1}
      \int_{x_n}^{x_{n+1}} \rho_{xx}\xi_{xx} \,dx +
      \sum_{n=0}^{N}\{\rho_{xxx}(x_n)\}[\xi(x_n)] -
      \sum_{n=0}^{N}\{\rho_{xx}(x_n)\}[\xi_x(x_n)].\label{Eq:n6}
    \end{equation}
    Now we estimate the first term as follows:
    \begin{eqnarray}
      \sum_{n=0}^{N-1} \int_{x_n}^{x_{n+1}} \rho_{xx}\xi_{xx}
      \,dx
      &\le& |||\rho_{xx}|||\,
            |||\xi_{xx}||| \le \frac{1}{6C}\|\rho\|_{H^2}^2 +
            \frac{C}{6} \|\xi\|_{\mathcal{E}}^2,\nonumber\\
      &\le& Ch^{2\min(k+1,s)-4}|||u|||_{H^s(\mathcal{E}_h)}^2 + \frac{C}{6}
            \|\xi\|_{\mathcal{E}}^2\label{Eq:6}.
    \end{eqnarray}
    Estimating the second term using H\"{o}lder's inequality, trace
    inequality and the Young's inequality, we obtain
    \begin{eqnarray}
      \sum_{n=0}^{N}\{\rho_{xxx}(x_n)\}[\xi(x_n)]
      \le Ch^{2\min(k+1,s)-6+\beta-1}|||u|||_{H^s(\mathcal{E}_h)}^2 + \frac{C}{6}
      \|\xi\|_{\mathcal{E}}^2\label{Eq:7}.
    \end{eqnarray}
    Similarly the last term can be estimated as
    \begin{eqnarray}
      \sum_{n=0}^{N}\{\rho_{xx}(x_n)\}[\xi_x(x_n)]
      &\le& Ch^{2\min(k+1,s)-6+\beta-1}|||u|||_{H^s(\mathcal{E}_h)}^2 + \frac{C}{6}
            \|\xi\|_{\mathcal{E}}^2\label{Eq:8}.
    \end{eqnarray}
    Combining \eqref{Eq:6}-\eqref{Eq:8}, we obtain
    the following bound for $\xi$ when  $\beta \ge 3$
    \begin{equation}
      \|\xi\|_{\mathcal{E}}\le Ch^{\min(k+1,s)-2} |||u|||_{H^s(\mathcal{E}_h)}.
    \end{equation}
    Now using $\|\eta\|_{\mathcal{E}} \le \|\xi\|_{\mathcal{E}} +
    \|\rho\|_{\mathcal{E}}$, we obtain the energy norm estimate
    for $\eta$. For the $L^2$-estimate of $\eta$, we use the Aubin
    Nitsch\'e duality argument. Consider the dual problem
    \begin{eqnarray*}
      &&\frac{d^4\phi}{dx^4} = \eta, \;\; x\in (a,b),\\
      &&\phi(a) = \phi(b) = 0,\\
      && \phi'(a) = \phi'(b) = 0. \qquad
    \end{eqnarray*}
    We note that $\phi$ satisfies the regularity condition $
    \|\phi\|_{H^4} \le C\|\eta\|$. Consider
    \begin{equation}
      (\eta,\eta) = (\eta,\phi_{xxxx}) = \sum_{n=0}^{N-1}
      \int_{x_n}^{x_{n+1}} \phi_{xx}\eta_{xx} \,dx +
      \sum_{n=0}^{N}\{\phi_{xxx}(x_n)\}[\eta(x_n)] -
      \sum_{n=0}^{N}\{\phi_{xx}(x_n)\}[\eta_x(x_n)].\nonumber
    \end{equation}
    Since $B(\eta,\chi) = 0 \;\; \forall \chi\in D^k(\mathcal{E}_h)$,
    we
    can write
    \begin{eqnarray}
      \|\eta\|^2 = (\eta,\eta) - B(\eta,\tilde{\phi})
      &=& \sum_{n=0}^{N-1}
          \int_{x_n}^{x_{n+1}} (\phi-\tilde{\phi})_{xx}\eta_{xx} \,dx +
          \sum_{n=0}^{N}\{(\phi-\tilde{\phi})_{xxx}(x_n)\}[\eta(x_n)]\nonumber\\
      &-&
          \sum_{n=0}^{N}\{(\phi-\tilde{\phi})_{xx}(x_n)\}[\eta_x(x_n)],\label{Eq:9}
    \end{eqnarray}
    where $\tilde{\phi}$ is a continuous interpolant of $\phi$ and
    satisfies the approximation property:
    \begin{equation}
      \|\phi-\tilde{\phi}\|_{H^q} \le
      Ch^{s-q}\|\phi\|_{H^s}.\label{eq:n7}
    \end{equation}
    We use the approximation property \eqref{eq:n7},
    the energy norm estimate for $\eta$ and the regularity result to
    bound each term on the right hand side of \eqref{Eq:9} and obtain
    the estimate
    for $\|\eta\|$ as:
    $$
    \|\eta\| \le Ch^{\min(k+1,s)}|||u|||_{H^s(\mathcal{E}_h)}.
    $$
    For the estimates of the temporal derivative of $\eta$, we
    differentiate \eqref{Eq:10} with respect to $t$ and repeat the
    arguments. Hence, it completes the rest of the proof.
  \end{proof}
\end{lemma}
\noindent
The following Lemma is useful to prove the error estimates:
\begin{lemma}\label{lemma1}
  Let $v\in \mathbb{P}_k(I_n)$ where $I_n = (x_n, x_{n+1}) \in
  \mathcal{E}_h$. Then there exists a positive constant $C$ independent of $h$
  such that,
  \begin{equation*}
    Ch_n^{-4}\|v\|_{L^2(I_n)}^2 \le \|v_{xx}\|^2_{L^2(I_n)},
  \end{equation*}
  where $h_n = \left(x_{n+1}-x_n\right)$.
  \begin{proof}
    We define the reference element $\hat{I_n}$ as
    \begin{equation*}
      \hat{I_n} = \left\{\left(\frac{1}{h_n}\right)x, \quad x
        \in I_n\right\}.
    \end{equation*}
    Since $v\in \mathbb{P}_k(I_n)$, we have the following relation
    (refer \cite{brenner}) for the norms in the reference element and the
    interval $I_n$
    \begin{equation*}
      |\hat{v}|_{H^r(\hat{I_n})} = h_n^{r - \frac{d}{2}}|v|_{H^r(I_n)}.
    \end{equation*}
    In one space dimension, i.e., $d=1$, we have
    \begin{eqnarray}
      \begin{aligned}
        \|\hat{v}\|_{L^2(\hat{I_n})} &=&
        h_n^{-\frac{1}{2}}\|v\|_{L^2(I_n)}, \;\; \mbox{and} \\
        |\hat{v}|_{H^2(\hat{I_n})} &=& h_n^{\frac{3}{2}}|v|_{H^2(I_n)}.
      \end{aligned}\label{neq:3.17}
    \end{eqnarray}
    By the equivalence of norm (refer \cite{brenner}), we have
    \begin{equation}
      C_0\|\hat{v}\|_{L^2(\hat{I_n})} \le \|\hat{v}_{xx}\|_{L^2(\hat{I_n})}
      \le C_1\|\hat{v}\|_{L^2(\hat{I_n})}.\label{neq:3.18}
    \end{equation}
    Now from \eqref{neq:3.17} and \eqref{neq:3.18}, we obtain
    \begin{equation*}
      C_0h_n^{-\frac{1}{2}}\|v\|_{L^2(I_n)} \le
      h_n^{\frac{3}{2}}\|v_{xx}\|_{L^2(I_n)}.
    \end{equation*}
    Rearranging the terms and squaring on both sides, we obtain the
    desired estimate.
  \end{proof}
\end{lemma}
\noindent
To obtain the error estimates, we subtract
\eqref{Neq2.1} from \eqref{Eq2.7} and using the auxiliary projection
\eqref{Eq:10}, we obtain the following error equation
\begin{equation}
  (\theta_t,\chi) + B(\theta_t,\chi) = (\eta_t,\chi) +
  (f(u^h)_x-f(u)_x,\chi).\label{Eq:12}
\end{equation}
Now we state and prove the following theorem.
\begin{theorem}
  Let $u^h(t)$ and $u(t)$ be the solutions of \eqref{Neq2.1} and
  \eqref{Eq2.7}, respectively. Let $u^h_0$ be the
  elliptic projection of $u_0$, i.e., $u^h_0 = \tilde{u}(0)$. Then for
  $s > 7/2$ and  
  there exists a positive constant $C$ independent of $h$ such that 
  \begin{eqnarray*}
    \|u(t)-u^h(t)\|_{\mathcal{E}} &\le& Ch^{\min(k+1,s)-2}
                                        \|u\|_{H^1(0,T;H^s(\mathcal{E}_h))},\\
    \|u(t)-u^h(t)\| &\le& Ch^{\min(k+1,s)}
                          \|u\|_{H^1(0,T;H^s(\mathcal{E}_h))}.
  \end{eqnarray*}
  \begin{proof}
    Setting $\chi = \theta(t)$ in \eqref{Eq:12}, we obtain
    \begin{equation}
      (\theta_t,\theta) + B(\theta_t,\theta) = (\eta_t,\theta) +
      (f(u^h)_x-f(u)_x,\theta).\label{Eqno:1}
    \end{equation}
    Now, we write equation \eqref{Eqno:1} as
    $$
    \frac{1}{2}\frac{d}{dt} \left( \|\theta(t)\|^2 +
    B(\theta(t),\theta(t)) \right) = (\eta_t,\theta) +
    (f(u^h)_x-f(u)_x,\theta).
    $$
    Integrating with respect to $t$ from $0$ to $T$ and noting that
    $\theta(0)=0$, we obtain
    \begin{equation}
      \|\theta(t)\|^2 + B(\theta(t),\theta(t)) = \int_{0}^{T} (\eta_t,\theta)\,ds
      + \int_0^T (f(u^h)_x-f(u)_x,\theta)\, ds.\label{Neq:8}
    \end{equation}
    We use integration by parts on the nonlinear term to obtain,
    \begin{equation}
      \left( (f(u^h) - f(u))_x,\theta\right) = \sum_{n=0}^N
      \{f(u^h)-f(u)\}[\theta(x_n)] + [f(u^h)-f(u)]\{\theta(x_n)\} +
      \left( f(u) - f(u^h),\theta_x\right).\label{Neq:5}
    \end{equation}
    Using the Cauchy Schwarz's and Young's inequality, we bound the last term of
    \eqref{Neq:5} as
    $$
    \left( f(u)-f(u^h),\theta_x \right) \le C(\|\eta\|^2 +
    \|\theta\|^2 + \|\theta\|_{\mathcal{E}}^2).
    $$
    Now for the first term in \eqref{Neq:5}, we use H\"{o}lder's
    inequality to write
    \begin{equation}
      \sum_{n=0}^N \{f(u^h)-f(u)\}[\theta(x_n)] \le
      \left(\sum_{n=0}^N|\{f(u^h)-f(u)\}|^2
      \right)^{1/2}\left(\sum_{n=0}^N
        |[\theta(x_n)]|^2\right)^{1/2}.\label{Neq:6}
    \end{equation}
    As earlier in \eqref{Eq:7}, we use the penalty term to write
    \eqref{Neq:6} as
    \begin{equation}
      \sum_{n=0}^N \{f(u^h)-f(u)\}[\theta(x_n)] \le
      Ch^{\frac{\beta-1}{2}}\|e\| \|\theta \|_{\mathcal{E}} \le
      C\|e\|\|\theta\|_{\mathcal{E}}, \quad \text{since} \;\; \beta
      \ge 3. \nonumber
    \end{equation}
    A similar bound for the second term can be obtained as
    follows. Using the H\"{o}lder's inequality, we write
    \begin{equation}
      \sum_{n=0}^N[f(u^h) - f(u)]\{\theta(x_n)\} \le
      \left(\sum_{n=0}^N h^{\beta}|[f(u^h)-f(u)]|^2
      \right)^{1/2}\left(\sum_{n=0}^N
        h^{-\beta}|\{\theta(x_n)\}|^2\right)^{1/2}.\label{Neq:7}
    \end{equation}
    Using the trace inequality, we obtain
    \begin{eqnarray*}
      \sum_{n=0}^N[f(u^h) - f(u)]\{\theta(x_n)\}
      &\le& \left(Ch^{\frac{\beta-1}{2}}\|e\|\right)
            \left(\sum_{n=0}^N
            h^{-\beta-1}\|\theta\|_{L^2(I_n)}^2\right)^{\frac{1}{2}}, \\
      &\le& Ch^{\frac{\beta-1}{2}}\left(\sum_{n=0}^N
            h^{-4 +
            (3-\beta)}\|\theta\|_{L^2(I_n)}^2\right)^{\frac{1}{2}}\|e\|,
            \\
      &\le& Ch \left(\sum_{n=0}^N
            h^{-4}\|\theta\|_{L^2(I_n)}^2\right)^{\frac{1}{2}}\|e\|
            \;\;\le\;\; C \|e\| \|\theta\|_{\mathcal{E}},
    \end{eqnarray*}
    where the last step is obtained by using \textsc{Lemma}
    \textbf{\ref{lemma1}}. Using the triangle
    inequality together with Young's inequality we obtain the bound
    \begin{eqnarray*}
      \left( (f(u^h) - f(u))_x,\theta\right)
      &\le& C(\|\theta\|^2 + \|\theta\|_{\mathcal{E}}^2) +
            Ch^{2\min(k+1,s)}|||u|||_{H^s(\mathcal{E}_h)}^2.
    \end{eqnarray*}
    Now using the coercivity of $B(.,.)$ and estimate of the
    nonlinear term in \eqref{Neq:8}, we arrive at
    \begin{eqnarray}
      \|\theta\|^2 + \|\theta\|_{\mathcal{E}}^2 \le
      Ch^{2\min(k+1,s)}\int_0^{T}|||u|||_{H^s(\mathcal{E}_h)}^2\, ds
      &+& Ch^{2\min(k+1,s)}\int_0^T|||u_t|||_{H^s(\mathcal{E}_h)}^2\, ds\nonumber\\
      &+& C\int_{0}^{t}(\|\theta\|^2 + \|\theta\|_{\mathcal{E}}^2) \,ds.
    \end{eqnarray}
    An application of Gronwall's inequality yields an estimate for
    $\theta$. We then use the triangle inequality to obtain the
    estimates for $e = u - u^h$. The estimates are optimal in
    $L^2$-norm if $\beta \ge 3$.
  \end{proof}
\end{theorem}

\section{Fully Discrete DGFEM}
\setcounter{equation}{0}


In this section, we derive a fully discrete DGFEM and establish
\textit{a priori} bounds
along with optimal error estimates. \\\\
{\bf Backward Euler discretization}:
Let $\Delta t$ denote the size of time discretization. Divide $[0,T]$
by \[ 0 = t_0 < t_1 < t_2 < \dots < t_{M-1} < t_M = T \] where
$t_{i+1} = t_i + \Delta t, \quad i = 0,1,\dots,M-1$ and $\Delta t =
\frac{T}{M}$. Let $u(t_n) = u^n$ and approximate $\frac{\partial
  u}{\partial t}$ by using backward Euler difference formula as
: \[\partial_t u^{n} = \frac{u^{n}-u^{n-1}}{\Delta t}.\] Now,
the fully discrete discontinuous Galerkin finite element method is
given as follows:
\begin{equation}
  \left(\partial_t U^{n+1}, \chi\right) + B(\partial_t U^{n+1},
  \chi) = \left(f(U^{n+1}), \chi\right), \;\; \mbox{for all} \;\; \chi \in
  D^k(\mathcal{E}_h), \label{Eqn3.1}
\end{equation}
\[U^0 = U_0, \] where $U^n$ is the fully discrete approximation of
$u(x,t_n)$.
\subsection{\textit{A priori} bounds}
In this sub-section, we prove an \textit{a priori} bound for the fully
discrete DGFEM.

\begin{theorem}
  Let $U^n$ be a solution to \eqref{Eqn3.1} and assume that $f'$ is
  bounded. Then there exists a positive constant $C$ such that
  \begin{equation}
    \|U^n\| + \|U^n\|_{\mathcal{E}} \le C(\|U_0\|).
  \end{equation}
  \begin{proof}
    Set $\chi = U^{n+1}$ in \eqref{Eqn3.1} to obtain
    \[\left(\frac{U^{n+1}-U^n}{\Delta t}, U^{n+1} \right) +
      B\left(\frac{U^{n+1}-U^{n}}{\Delta t}, U^{n+1}\right) = \left(
        f(U^{n+1})_x, U^{n+1} \right).\] Multiplying by $\Delta
    t$ throughout, we arrive at
    \begin{equation}
      \left(U^{n+1}-U^n, U^{n+1} \right)  +
      B\left(U^{n+1}-U^{n}, U^{n+1}\right) = \Delta t \left(
        f(U^{n+1})_x, U^{n+1} \right). \label{Eqn3.2}
    \end{equation}
    Using the fact that for any two real numbers $x$ and $y$, we have
    \[
      \frac{1}{2}(x^2-y^2) \le \frac{1}{2} (x^2-y^2+(x-y)^2) = (x-y)x. \]
    We rewrite the equation (\ref{Eqn3.2}), we obtain
    \begin{equation}
      \frac{1}{2}\lVert U^{n+1}\rVert^2 - \frac{1}{2}\lVert
      U^{n} \rVert^2 +
      \frac{1}{2}\left(B\left(U^{n+1}, U^{n+1}\right) - B(U^n,U^n)\right) \le \Delta t \left(
        f(U^{n+1})_x, U^{n+1} \right).\label{eq:3-40}
    \end{equation}
    \noindent
    Using the Cauchy Schwarz's inequality, Poincare inequality and Young's
    inequality with the bound on $f'(u)$, we obtain the
    following inequality from \eqref{eq:3-40}
    \begin{equation*}
      \frac{1}{2}\lVert U^{n+1}\rVert^2 - \frac{1}{2}\lVert
      U^{n} \rVert^2 + \frac{1}{2}\left( B(U^{n+1},U^{n+1}) - B(U^n,U^n)\right) \le C \Delta t \left(
        \frac{\lVert U^{n+1}\rVert_{\mathcal{E}}^2}{2} +\frac{\lVert
          U^{n+1} \rVert^2}{2} \right).
    \end{equation*}
    Summing over $n = 0,1,2,\dots,J-1$ and using the coercivity of $B(\cdot,\cdot)$, we obtain
    \begin{equation*}
      \frac{1}{2}\lVert U^{J}\rVert^2 + \frac{1}{2}\lVert
        U^{J} \rVert_{\mathcal{E}}^2 \le \frac{1}{2}\lVert U^{J}\rVert^2 + \frac{1}{2}B(U^0,U^0) +  C \Delta t \sum_{n=0}^{J-1}\left(
        \lVert U^{n+1}\rVert_{\mathcal{E}}^2 +\lVert U^{n+1} \rVert^2 \right).
    \end{equation*}
    Rearranging the terms and applying the discrete Gronwall Inequality,
    we obtain the desired \textit{a priori} bound on $U^J$.
  \end{proof}
\end{theorem}

\subsection{Error Estimates}

In this sub-section, we prove the optimal error estimates for the fully discrete
discontinuous Galerkin method. \\\\
Subtracting equation \eqref{Eqn3.1} from \eqref{Eq2.7} and using the
auxiliary projection \eqref{Eq:10}, we obtain the error equation as
\begin{eqnarray}
  \left(\partial_t \theta^{n+1}, \chi\right) + B(\partial_t
  \theta^{n+1}, \chi)
  &=& \left(\partial_t \eta^{n+1}, \chi\right) +
      \left((f(U^{n+1})_x - f(u^{n+1})_x),\chi\right)\nonumber\\
  &+&
      \left(\sigma^{n+1},\chi\right) + B(\sigma^{n+1},\chi). \label{eq:3-45}
\end{eqnarray}
where $\sigma^{n+1} = u_t^{n+1}-\partial_t U^{n+1}$. Before we derive
the error estimate, we state and prove the following Lemma.
\begin{lemma}\label{Lemma2}
  Let $\sigma^{n} = u_t^{n} - \partial_t U^{n}$. Then the following
  holds
  \[ \lVert \sigma^n \rVert^2 \le \Delta t \int_{t_{n-1}}^{t_n} \lVert
    u_{tt}(s)\rVert^2 ds.\]
  \begin{proof}
    Consider $I = \int_{t_{n-1}}^{t_{n}} (s-t_{n-1}) u_{tt}(s) \:
      ds$. Using integration by parts, we see that $I = \Delta t
      \; \sigma^{n}$. A use of the H\"{o}lder's inequality, we obtain
    $\lVert \sigma^n \rVert^2 \le \Delta t \int_{t_{n-1}}^{t_n} \lVert
      u_{tt}(s)\rVert^2 ds $.
    The same inequality can be proved in the energy norm as well.
  \end{proof}
\end{lemma}

\begin{theorem}
  Let $U^0 = \tilde{u}(0)$ so that $\theta^0 = 0$. Then there exists a
  positive constant $C$ independent of $h$ and $\Delta t$ such that,
  \begin{eqnarray}
    \lVert u(t_n) - U^n \rVert_{\mathcal{E}}
    &\le& C\left(h^{\min(k+1,s)-2}\,\|u\|_{H^1(0,T,H^s(\mathcal{E}_h))}+
          \Delta t\,
          |||u_{tt}|||_{L^2(0,T;H^s(\mathcal{E}_h))}\right),\nonumber\\
    \lVert u(t_n) - U^n \rVert
    &\le& C\left(h^{\min(k+1,s)}\,\|u\|_{H^1(0,T,H^s(\mathcal{E}_h))} +
          \Delta t\, |||u_{tt}|||_{L^2(0,T;H^s(\mathcal{E}_h))}\right).
  \end{eqnarray}
  \begin{proof}
    Set $\chi = \theta^{n+1}$ in \eqref{eq:3-45} we obtain
    \begin{eqnarray}
      \left(\theta^{n+1}-\theta^n, \theta^{n+1}\right) + B(\theta^{n+1}-\theta^n,
      \theta^{n+1})
      &=& \left(\eta^{n+1}-\eta^n, \theta^{n+1}\right) +
          \Delta t \left(f(U^{n+1})_x - f(u^{n+1})_x,\theta^{n+1}\right)\nonumber\\
      &+&
          \Delta t \left(\sigma^{n+1},\theta^{n+1}\right)  + \Delta  t
          B(\sigma^{n+1},\theta^{n+1}).\label{eq:3-46}
    \end{eqnarray}
    Using Cauchy Schwarz's Inequality and
    constructing upper bounds similar to the semidiscrete case on the
    right hand side of \eqref{eq:3-46}, we obtain the inequality
    \begin{eqnarray}
      \frac{1}{2}\lVert \theta^{n+1} \rVert^2
      -
      \frac{1}{2}\lVert \theta^{n} \rVert^2 + \frac{1}{2}\left(
      B(\theta^{n+1},\theta^{n+1}) - B(\theta^n,\theta^n) \right)
      \le  Ch^{2\min(k+1,s)}\int_{t_n}^{t_{n+1}}\|u_t\|^2\, ds \nonumber\\
      +\, Ch^{2\min(k+1,s)+\beta-3}|||u|||_{H^s(\mathcal{E}_h)}^2
      + C\Delta t \left( \lVert \theta^{n+1}
      \rVert^2 + \lVert \theta^{n+1}
      \rVert_{\mathcal{E}}^2 +
      \lVert \sigma^{n+1}
      \rVert^2 + \lVert \sigma^{n+1}
      \rVert_{\mathcal{E}}^2\right).\label{eqNo:11}
    \end{eqnarray}
    Sum over $n=0,1,2,\dots,M-1$ on both sides of \eqref{eqNo:11} and
    using \textsc{Lemma} \textbf{\ref{Lemma2}}, we obtain
    \begin{eqnarray}
      \lVert \theta^{M} \rVert^2 + \lVert \theta^{M}
      \rVert_{\mathcal{E}}^2
      &\le& Ch^{2\min(k+1,s)}|||u|||_{H^s(\mathcal{E}_h)}^2 +
            Ch^{2\min(k+1,s)}\int_{0}^{T}\|u_t\|^2\, ds
            \nonumber \\
      &+& C\Delta t
          \sum_{n=0}^{M-1}\left( \lVert
          \theta^{n+1}
          \rVert + \lVert \theta^{n+1}
          \rVert_{\mathcal{E}}^2\right) + C \Delta
          t^2\|u_{tt}\|^2_{L^2(0,T;H^s(\mathcal{E}_h))}.
    \end{eqnarray}
    Using the discrete Gronwall's inequality, we obtain the estimate
    for $\theta^M$ as
    \[ \lVert \theta^M \rVert^2 + \lVert \theta^{M}
      \rVert_{\mathcal{E}}^2 \le
      C\left(h^{2\min(k+1,s)}\|u\|_{H^1(0,T,H^s(\mathcal{E}_h))}^2 + \Delta
        t^2\|u_{tt}\|^2_{L^2(0,T;H^s(\mathcal{E}_h))}\right),\]
    provided $\beta \ge 3$.
    Using the triangle inequality, we can prove the required estimate for
    $\lVert e \rVert_{\mathcal{E}}$ and $\lVert e \rVert$.
  \end{proof}
\end{theorem}
\section{Numerical Results}
\setcounter{equation}{0}

\begin{figure}[H]
  \begin{subfigure}[b]{0.5\textwidth}
    \includegraphics[height=3.25cm,width=7.6cm]{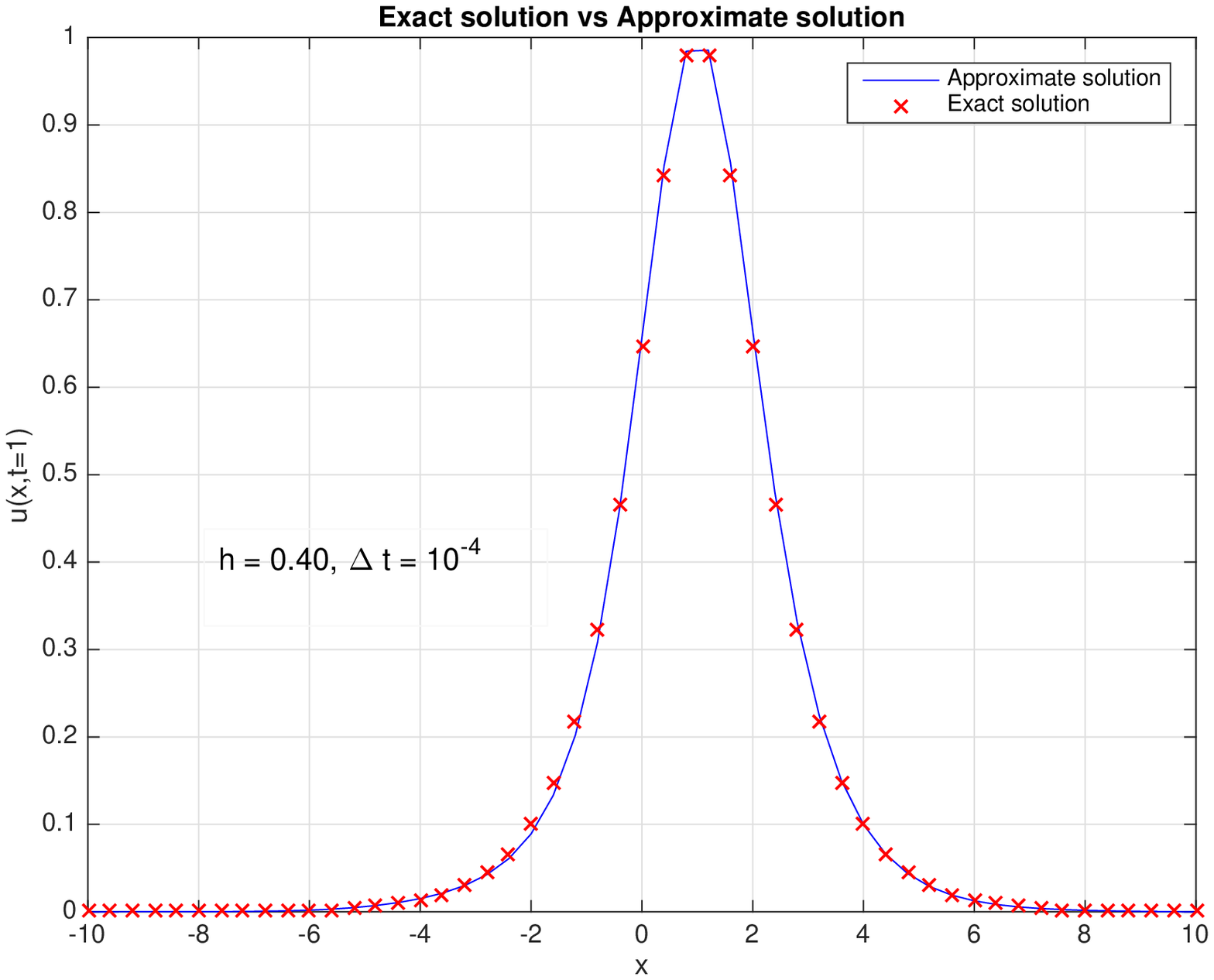}
    \caption{$N=50$}
    \label{fig:1}
  \end{subfigure}
  \begin{subfigure}[b]{0.5\textwidth}
    \includegraphics[height=3.25cm,width=7.6cm]{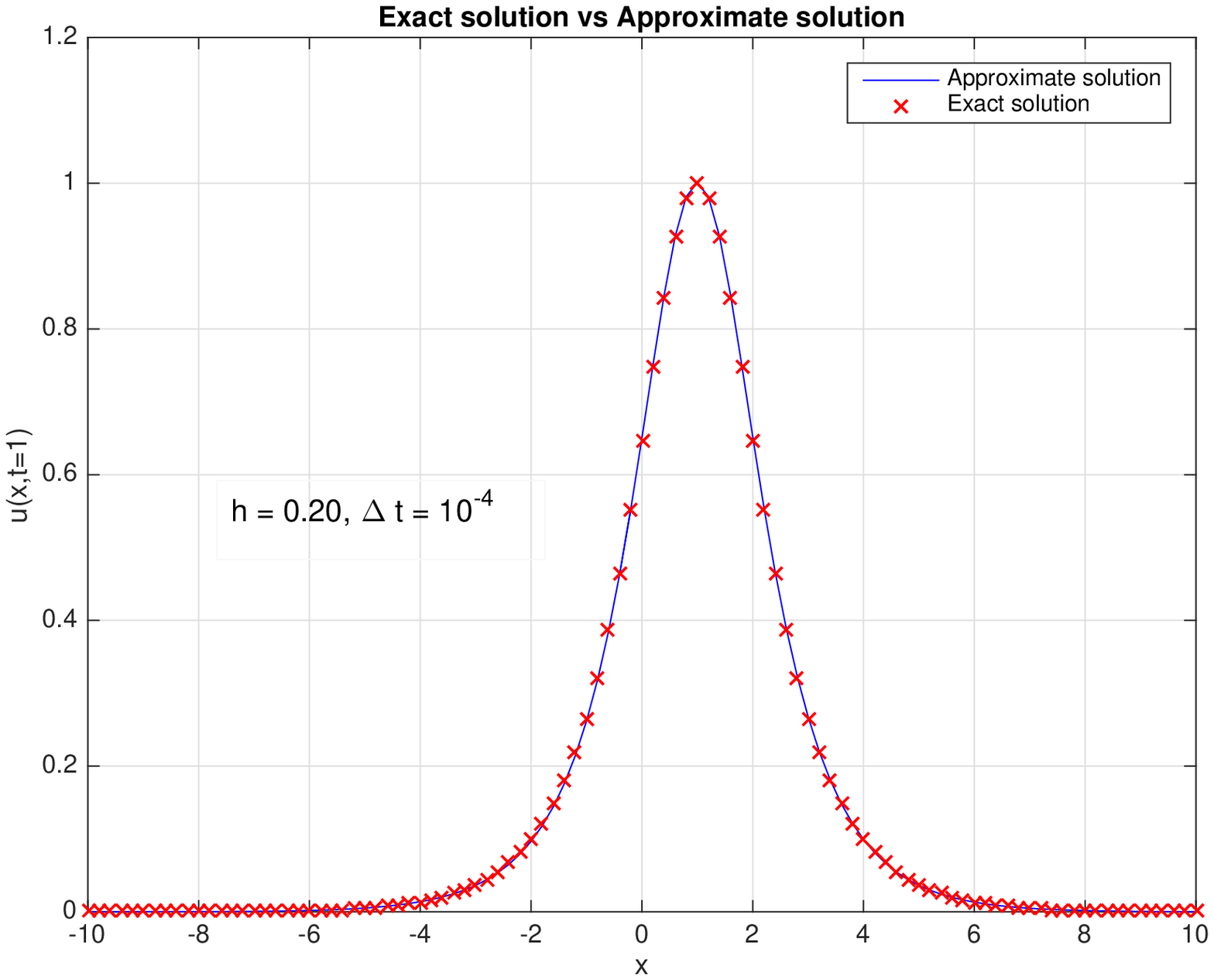}
    \caption{$N=100$}
    \label{fig:2}
  \end{subfigure}
  \begin{subfigure}[b]{0.5\textwidth}
    \includegraphics[height=3.25cm,width=7.6cm]{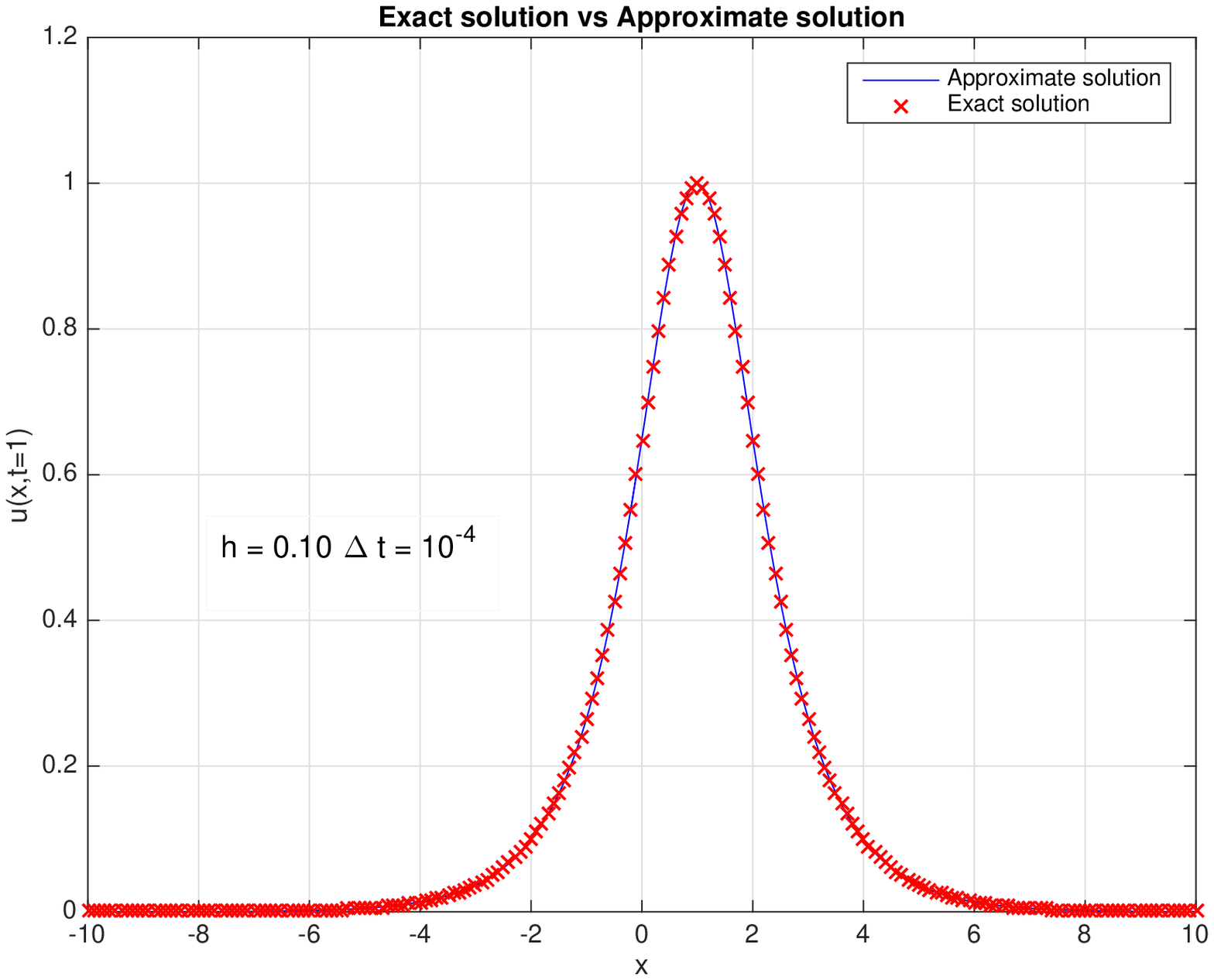}
    \caption{$N=200$}
    \label{fig:3}
  \end{subfigure}
  \begin{subfigure}[b]{0.5\textwidth}
    \includegraphics[height=3.25cm,width=7.6cm]{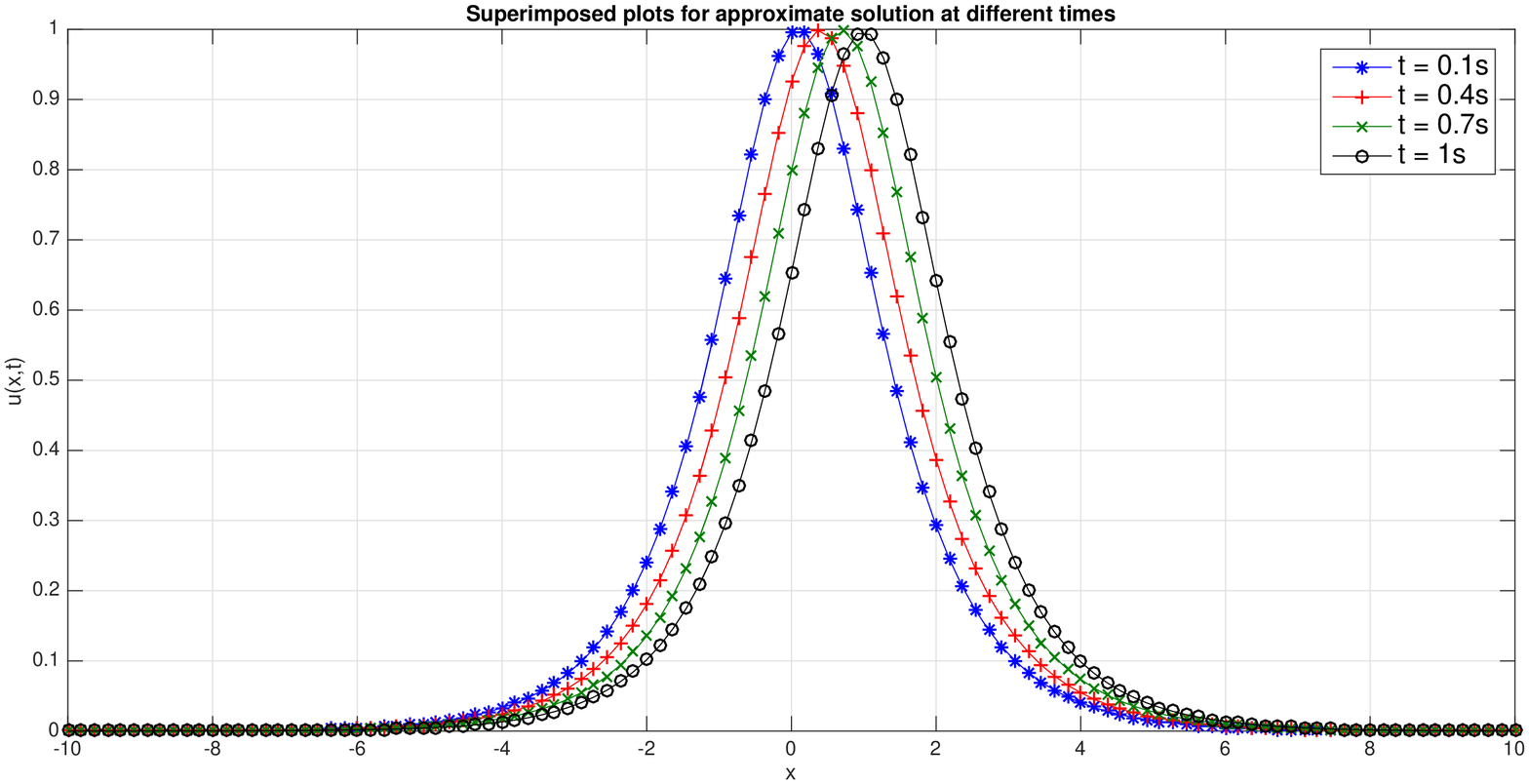}
    \caption{Solution profiles from $t=0$ to $t=1$.}
    \label{fig4}
  \end{subfigure}
  \caption{Approximate solution using $P_2$ elements with exact solution $u(x) = \text{sech}(x-1)$}
\end{figure}

In this section, we perform some numerical experiments to validate the
theoretical results.  \\\\
We consider the following Rosenau equation
\begin{equation}
  u_t + \frac{1}{2}u_{xxxxt} = f(u)_x \label{Eq5.1}
\end{equation}
with the boundary conditions (\ref{Eq1.3}) and the nonlinear function $f(u) =
10u^3-12u^5-\frac{3}{2}u$. The exact solution of the 
equation \eqref{Eq5.1} is $u(x,t) = \text{sech}(x-t)$. Since the
equation \eqref{Eq5.1} was considered as a benchmark example to
validate the results by several authors (for instance
\cite{manickam,dgros}), the same example has been taken to compare the
existing results.\\\\
We choose the computational domain $\Omega = (-10,10)$ and the final
time $T = 1$. The equation is solved numerically with corresponding initial and
boundary conditions.\\\\
The order of convergence for the numerical method was computed by
using the formula
\begin{equation*}
  p \approx \frac{\log\left(\frac{\lVert E_i \rVert}{\lVert E_{i + 1}
        \rVert}\right)}{\log\left(\frac{h_i}{h_{i + 1}}\right)}, \; i
  = 1, 2, 3, 4.
\end{equation*}
In Table \ref{tab:1}, we show the order of convergences for piecewise
quadratic and piecewise cubic basis functions.
\begin{table}[H]
  \begin{center}
    \begin{tabular}{|c|c|c|c|c|c|}
      \hline
      \multicolumn{3}{|c|}{Quadratic Elements ($k=2$)} &
                                                         \multicolumn{3}{|c|}{Cubic Elements ($k=3$)}\\
      \hline
      $h$ & Error $\lVert e \rVert$ & Order & $h$ &Error $\lVert e \rVert$ & Order \\
      \hline
      $0.20000$ & $1.47828\times 10^{-2}$ &- & $0.40000$& $5.17824\times 10^{-2}$ & -\\
      \hline
      $0.18182$ & $1.12327\times 10^{-2}$ & $2.8815$ & $0.33331$ &$2.58529\times10^{-2}$ &$3.8099$\\
      \hline
      $0.16667$ & $8.68469\times 10^{-3}$ & $2.9567$ & $0.28571$ &$1.41359\times10^{-2}$ &$3.9164$\\
      \hline
      $0.15384$ & $6.87521\times 10^{-3}$ & $2.9189$ & $0.25000$ &$8.35109\times10^{-3}$ &$3.9416$\\
      \hline
      $0.14257$ & $5.46407\times 10^{-3}$ & $3.0999$ & $0.22222$ & $5.23371\times10^{-3}$ & $3.9672$\\
      \hline
    \end{tabular}
  \end{center}
  \caption{Order of convergence for $P_2$ and $P_3$ elements with
    $\sigma_0 = \sigma_1 = 2000$.}
  \label{tab:1}
\end{table}
\noindent
Below, we show the Figures \eqref{fig:1}-\eqref{fig:3}, for the
comparison of exact solution profile with that of the approximate
solution obtained from DGFEM. The Figure \eqref{fig4} shows the solution
profile at different time levels.
\noindent
We compare our numerical results with Choo {\it et al.}
\cite{dgros}. We observe that our solution profiles matches very
acurately and we have achieved third order convergence for quadratic
elements and fourth order convergence for cubic elements  which are
optimal. The proposed method can be easily
  extended to higher degree polynomials and higher dimensions also. 
  But the cGdG method
  considered in \cite{dgros} is difficult to apply for higher dimensions due to the requirement of $C^1$-elements.
\subsection{Decay Estimates}
\begin{figure}[H]
  \centering
  \includegraphics[height=8cm,width=14cm]{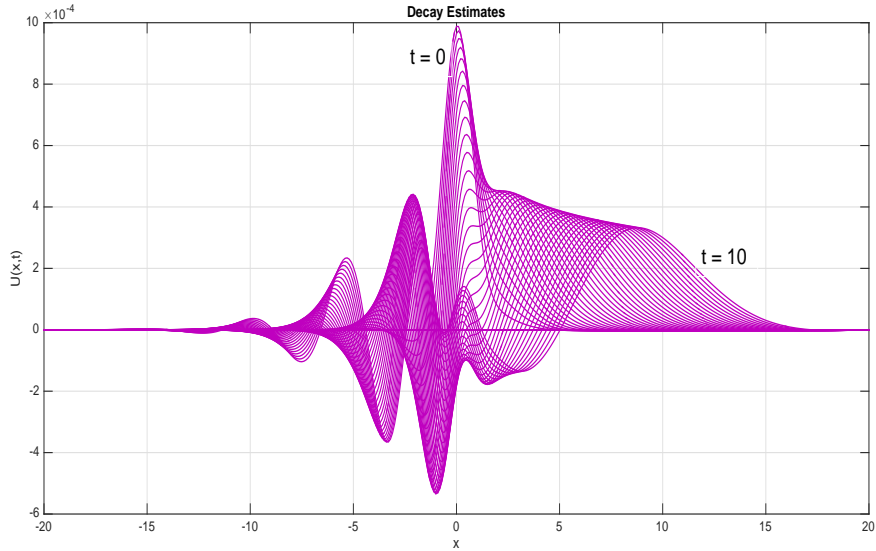}
  \caption{Curves illustrating the decaying nature of the
    solution}
  \label{fig:4}
\end{figure}

\begin{figure}[H]
  \centering
  \includegraphics[height=8cm,width=14cm]{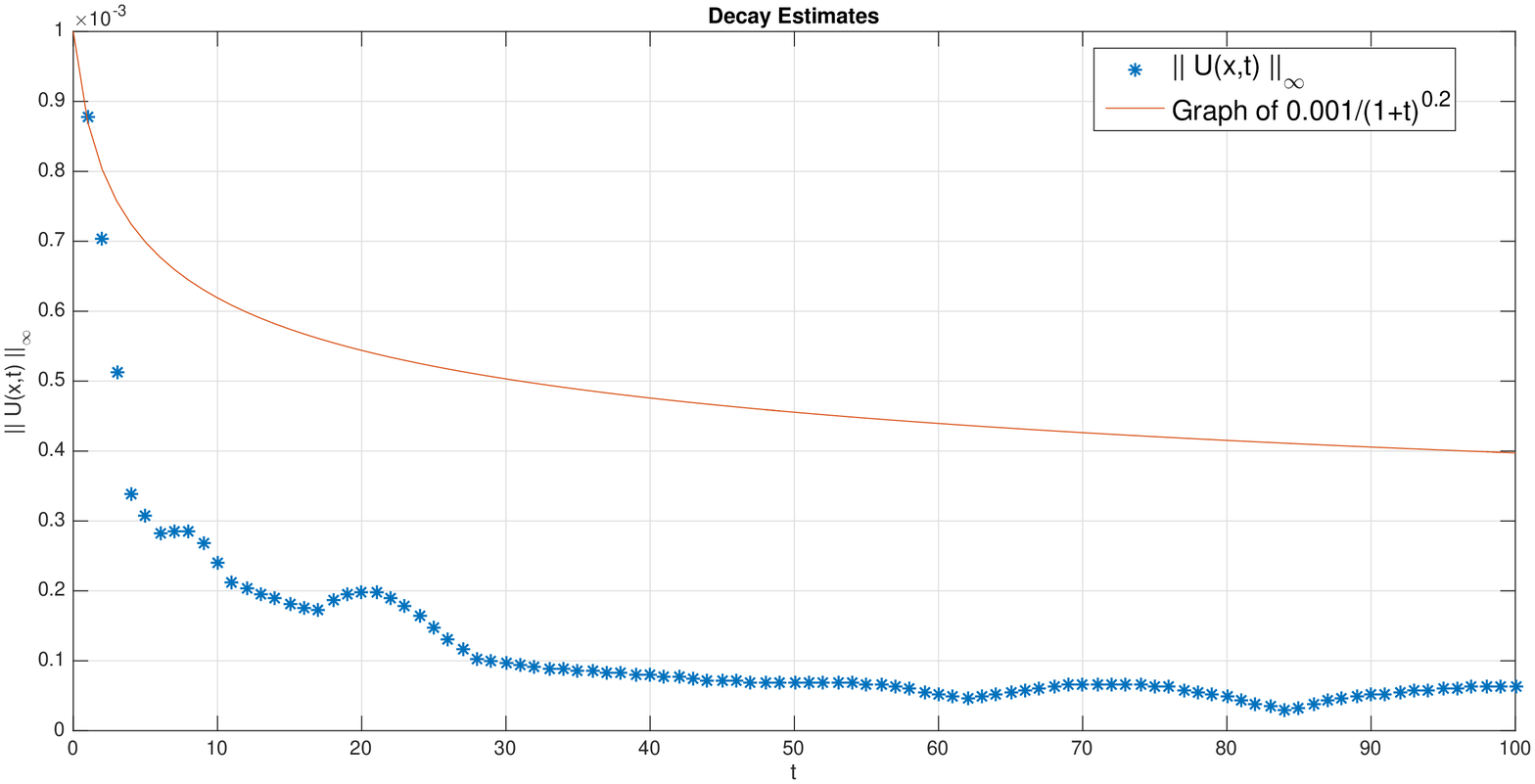}
  \caption{The variation of the  $L^\infty$-norm of the
    solution with respect to time ($t=0-100$).}
  \label{fig:5}
\end{figure}

In this sub-section, we validate the decay estimates that was derived by
Park in \cite{park2}. As in \cite{manickam}, we consider the following equation
\begin{eqnarray}
  u_t + u_{xxxxt} + u_x = f(u)_x, \quad (x,t) \in (0,1) \times
  (0,T]\label{eq:rosenau}
\end{eqnarray}
with the initial condition
\begin{eqnarray}
  u(x,0) = \phi_0(x), \label{eq:init}
\end{eqnarray}
and boundary conditions (\ref{Eq1.3}). Here $f(u) = \sum_{i=1}^n \frac{c_iu^{p_i+1}}{p_i+1}, \; c_i \in
  \mathbb{R}, p_i>0.$ \\\\
For a small initial data, it has been proven that the solution to the
Rosenau equation with small initial data decays like $\frac{1}{(1+t)^{\frac{1}{5}}}$ in the
$L^\infty$-norm for $\underset{1\le i\le n}{\min} p_i\; > 6$. Like in \cite{manickam}, we take
$c_1 = -1, p_1 = 7$, $c_2 = \frac{4}{7}, p_2 = 8$, $c_3 = -\frac{4}{3}, p_3 = 9$ and
$\phi_0(x) = 0.001e^{-x^2}$ in \eqref{eq:init}. The solution curves for
\eqref{eq:rosenau} and \eqref{eq:init} for $t=0$ to $t=10$ is shown in
Figure \eqref{fig:4}. The height of the initial pulse decreases with
time, indicating a decaying behavior of the solution. The decay in
$L^\infty$-norm of the approximate solution with time is shown in
Figure \eqref{fig:5}.

\section{Conclusion}
In this paper, we derived \textit{a priori} bounds and optimal error
estimates for the semidiscrete problem. Next, we discretized the
semidiscrete problem in the temporal direction using a backward Euler
method, and derived \textit{a priori} bounds and optimal error
estimates. We have validated the theoretical results by performing
some numerical experiments. Compared to the existing results, our
method requires less regularity of the original problem.

\section{Acknowledgements}
The authors would like to thank Department of Science and Technology
DST-FIST Level-1 Program Grant No. SR/FST/MSI-092/2013 for providing
computational facilities.



\end{document}